\documentclass[12pt,twoside,doublespace]{article}
\usepackage{latexsym,amssymb}

\def\smskip{\par\vskip 5 pt}
\def\QED{\hfill $\Box$\smskip}

\newtheorem{theorem}{Theorem}[section]

\newtheorem{proposition}{Proposition}[section]

\newtheorem{remark}{Remark}[section]

\begin{document}

\begin{center}

\vspace{15pt}

{\Large \bf Application of Market  Models}

\vspace{10pt}

{\Large \bf  to Network Equilibrium Problems}

\vspace{35pt}

{\sc I.V.~Konnov\footnote{\normalsize Department of System Analysis
and Information Technologies, Kazan Federal University, ul.
Kremlevskaya, 18, Kazan 420008, Russia. E-mail: konn-igor@ya.ru}}

\end{center}

\begin{abstract}
We present a general two-side market model with divisible
commodities and price functions of participants.
A general existence result on unbounded sets is obtained
from its variational inequality re-formulation.
We describe an extension of the network flow
equilibrium problem with elastic demands and 
 a new equilibrium type model for resource
allocation problems in wireless communication networks, which appear to be
particular cases of the general market
model. This enables us to obtain new existence results for these models as some
adjustments of that for the market model.
Under certain additional
conditions the general market model can be reduced to a decomposable optimization problem
where the goal function is the sum of two  functions and one of them is
convex separable, whereas the feasible set is the corresponding
Cartesian product. We discuss some versions of the partial linearization method, which
can be applied to these network equilibrium problems.


\end{abstract}


\section{Introduction}\label{sc:1}

Investigation of complex systems with active elements having their
own interests and sets of actions is usually based on a
suitable equilibrium concept. Such a concept should equilibrate
different interests and opportunities of the elements (agents,
participants) and provide ways of its proper implementation within
some accepted basic (information) behavior framework of the system
under investigation.

For instance, the classical perfectly
(Walrasian) and imperfectly (Cournot - Bert\-rand) competitive
models, which are most popular in economics (see e.g.
\cite{Nik68,OS90} and references therein), reflect different
equilibration mechanisms and information frameworks.
We recall that actions of any separate agent within
a perfect competition model can not impact the state of the whole
system, hence any agent may utilize some integral system
parameters (say, prices), rather than
the information about the behavior of other agents.
On the contrary, actions of any separate
agent in an imperfectly competitive model can change the state of
the whole system. Therefore, the model is formulated as
a (non-cooperative) game problem and is usually based on
the well-known Nash equilibrium concept \cite{Nas51}.
Nevertheless, real systems (markets) may give wide variety of these
features and different information frameworks. Hence,
flexible equilibrium models could also be very useful for derivation of
efficient decisions in complex systems.

In this paper, we consider a general two-side market model with divisible
commodities and price functions of participants. It is based on the
  auction market models proposed in \cite{Kon06,Kon07a}, where
the equivalence result with a variational inequality problem was established.
Afterwards, some extensions to the
multi-commodity case and applications to resource allocation in
telecommunication networks were suggested in \cite{Kon07b,Kon13}.
The alternative equilibrium concept related to this model
was proposed in \cite{Kon15e}, where it was also shown
 that the same equilibrium state can be attained within
different mechanisms and information exchange schemes, including the
completely decentralized competitive mechanism.

We now suggest a somewhat more general class of market equilibrium models,
which follows the approach from \cite{Kon15e}. It is subordinated
to the material balance condition and can also  be formulated as a variational inequality problem,
hence one can utilize the well-developed theory and methods of variational inequalities
for investigation and solution finding of this equilibrium model.
We give a new existence result for the model in the case where the feasible set is
unbounded. Besides, under certain integrability
conditions the model can be also reduced to an optimization problem.
We suggest  a new cyclic version of the partial linearization method
for its decomposable case. We describe extensions of the known network flow
equilibrium problem with elastic demands and a resource
allocation problem in wireless communication networks and show they are
particular cases of the presented market
model.  This enables us to obtain new existence results for these models
 and to solve these problems with the partial linearization method.


\section{A general multi-commodity market equilibrium model}\label{sc:2}

We start our considerations from a general
market model with $n$ divisible commodities, which somewhat extends those in
\cite{Kon07b,Kon15e}; see also \cite{Kon16a} for the vector model.
 For each commodity $s$, each trader $i$ chooses
some offer value $x_{is}$ in his/her capacity segment $[\alpha'
_{is}, \alpha'' _{is}]$ and has a price function $g_{is}$.
Similarly, each buyer $j$ chooses some bid
value $y_{js}$ in his/her capacity segment $[\beta' _{js}, \beta
''_{js}]$ and has a price function $h_{js}$.
We denote by $I_{s}$ and $J_{s}$ the finite index sets of traders and buyers
attributed to commodity $s$ and set $N=\{1,\dots, n\}$.
Clearly, each trader/buyer can be attributed to many commodities.  We suppose that the
prices may in principle depend on all the bid/offer volumes of all the commodities. That
is, if we set $x_{(s)}=(x_{is})_{i \in I_{s}}$, $x=(x_{(s)})_{s \in N}$,
$y_{(s)}=(y_{js})_{j \in J_{s}}$, $y=(y_{(s)})_{s \in N} $, and
$w=(x,y)$,  then $g_{is}=g_{is}(w)$ and $h_{js}=h_{js}(w)$.
Let $b_{s}$ denote the value of the external excess demand
for  commodity $s$, then  $b=(b_{s})_{s \in N}$. If it equals zero, the market is closed.
Any market solution must satisfy the balance equation, hence we
obtain the feasible set of offer/bid values
\begin{eqnarray*}
   W &=& \prod_{s \in N} W_{s}, \ \mbox{where} \\
   W_{s}&=& \left\{ w_{(s)}=(x_{(s)},y_{(s)}) \ \vrule \
\begin{array}{l}
\sum \limits_{i \in I_{s}} x_{is} - \sum \limits_{j \in J_{s}} y_{js} = b_{s}; \\
x_{is} \in [\alpha' _{is}, \alpha'' _{is}], i \in I_{s}, y_{js} \in [\beta'
_{js}, \beta ''_{js}], j \in J_{s}
\end{array}
\right\}; \\
 && \mbox{for} \ s \in N.
\end{eqnarray*}
 A vector $\bar w= (\bar x, \bar y)\in W$ is said to be a {\em market
equilibrium point} if there exists a price vector $\bar p=(\bar p_{s})_{s \in N}$ such that
\begin{equation} \label{eq:2.1}
\begin{array}{c}
\displaystyle g_{is}(\bar w) \left\{
\begin{array}{ll}
\geq
\bar p _{s} & \quad \mbox{if} \quad \bar x_{is}=\alpha '_{is}, \\
=\bar p _{s} & \quad \mbox{if} \quad \bar x_{is} \in (\alpha
'_{is},
\alpha ''_{is}), \\
\leq \bar p _{s} & \quad \mbox{if} \quad \bar x_{is}=\alpha
''_{is},
\end{array}
 \quad \mbox{for} \ i \in I_{s}; \right.
\end{array}
\end{equation}
and
\begin{equation} \label{eq:2.2}
\begin{array}{c}
\displaystyle h_{js}( \bar w) \left\{
\begin{array}{ll}
\leq \bar p _{s} & \quad \mbox{if} \quad \bar y_{js}=\beta '_{js}, \\
=\bar p _{s} & \quad \mbox{if} \quad \bar y_{js} \in (\beta '_{js},
     \beta ''_{js}), \\
\geq \bar p _{s} & \quad \mbox{if} \quad \bar y_{js}=\beta ''_{js},
\end{array} \quad \mbox{for} \ j \in J_{s}; \right.
\end{array}
\end{equation}
for $s \in N$. We now give the basic relation between
the market equilibrium problem (\ref{eq:2.1})--(\ref{eq:2.2}) and a
variational inequality (VI, for short). Its proof is almost the same as that in
\cite[Theorem 2.1]{Kon07b} and is omitted.

\begin{proposition} \label{pro:2.1}

(a) If $(\bar w, \bar p)$ satisfies (\ref{eq:2.1})--(\ref{eq:2.2})
for $s \in N$ and $\bar w \in W$, then $\bar w$ solves VI:
Find $\bar w \in W$ such that
\begin{equation}\label{eq:2.3}
\sum \limits_{s \in N}  \left[ \sum \limits_{i \in I_{s}} g_{is} (\bar
w) (x_{is} - \bar x_{is}) - \sum \limits_{j \in J_{s}} h_{js} (\bar w)
(y_{js} - \bar y_{js})\right] \geq 0 \quad \forall w \in W.
\end{equation}

(b) If a vector $\bar w$ solves VI (\ref{eq:2.3}), then there exists
$\bar p \in \mathbb{R}^{n}$ such that $(\bar w, \bar p)$ satisfies
(\ref{eq:2.1})--(\ref{eq:2.2}) for $s \in N$.
\end{proposition}
The presence of the price functions is invoked by
complexity of the whole system, i.e., the price functions
may contain participants' intentions or
reflect interdependence (mutual influence) of the elements,
which need not be known to the participants. 

It follows from Proposition \ref{pro:2.1} that we can establish
existence results for equilibrium problems of form
(\ref{eq:2.1})--(\ref{eq:2.2}) by using suitable results from the
theory of VIs or more general equilibrium problems.
For instance, if the feasible set $W$ is bounded
and  the cost mapping of VI (\ref{eq:2.3}) is continuous,
then equilibrium problem (\ref{eq:2.1})--(\ref{eq:2.2}) has a solution.
In the unbounded case, we need certain coercivity assumptions.

We follow the approach from \cite{Kon15e,Kon16a} and consider for simplicity the
case where all the lower bounds $\alpha '_{is}$ and $\beta '_{js}$ of
capacities are fixed and greater than $-\infty$, whereas some
upper bounds $\alpha ''_{is}$ and $\beta ''_{js}$ can be absent.
Then, for each commodity $s \in N$, we define the index sets
$$
I^{u}_{s}=\{ i \in I_{s} \ | \ \alpha ''_{is}=+\infty \} \ \mbox{and} \
J^{u}_{s}=\{ j \in J_{s} \ | \ \beta ''_{js}=+\infty \},
$$
and take the following coercivity condition.

{\bf (C)} {\em There exists a number $r>0$ such that for any point
$w=(x,y) \in W$ and for each $s \in N$ it holds that}
\begin{eqnarray*}
 && \forall l \in J^{u}_{s},  y_{ls}>\max \{r, \beta '_{js}\} \Longrightarrow  \exists k \in I^{u}_{s} \  \mbox{such that} \\
  &&  x_{ks}>\alpha '_{ks} \ \mbox{and} \  g_{ks}(w)\geq h_{ls}(w).
\end{eqnarray*}
This condition seems rather natural: at any
feasible point $w$ and for each fixed commodity $s$, any large
 demand value of buyer $l$ invokes existence of a trader
$k$ whose price is not less than the price  of buyer $l$.

\begin{proposition} \label{pro:2.2} Suppose that the set $W$
is nonempty, all the functions $g_{is}$ and $g_{js}$ are continuous
for all $i \in I$, $j \in J$, and $s \in N$. If condition {\bf (C)}
is fulfilled, then VI (\ref{eq:2.3}) has a solution.
\end{proposition}
The proof of this assertion is almost the same as those in
\cite[Theorem 1]{Kon15e} and \cite[Theorem 4.3]{Kon16a} and is omitted.


\section{Partial linearization methods} \label{sc:3}

Due to Proposition \ref{pro:2.1}, we can take
various iterative solution methods for optimization and
variational inequality problems (see e.g. \cite{Kon07a,Kon15e,Kon13a})
for finding solutions of the market
equilibrium problems of form (\ref{eq:2.1})--(\ref{eq:2.2}).
We now intend to consider a special integrable class of these problems
that admits efficient iterative solution methods.

Let us first take a problem of minimization of the sum of two functions $\mu(w)+\eta(w)$
over a feasible set $W \subseteq \mathbb{R}^{m}$, or briefly,
\begin{equation} \label{eq:3.1}
 \min \limits _{w \in W} \to \left\{\mu(w)+\eta(w)\right\}.
\end{equation}
We suppose that the set $W \subset \mathbb{R}^{m}$ is non-empty, convex, and compact,
both the functions are convex and $\mu : \mathbb{R}^{m} \to \mathbb{R}$ is smooth.
Moreover, the minimization of the function $\eta$ over the set $W$
is not supposed to be difficult. In this case one can apply the 
partial linearization (PL for short) method,
which was first proposed in \cite{MF81}.

\medskip \noindent
 {\bf Method (PL).} \\
 Choose a point $w^{0}\in W$ and set $k=0$. At the $k$-th iteration, $k=0,1,\ldots$, we have
a point $w^{k}\in W$. Find some solution $v^{k}$ of the problem
\begin{equation} \label{eq:3.2}
\min_{v \in W} \to \left\{\langle \mu'(w^{k}),v\rangle + \eta(v)\right\}
\end{equation}
and define $p^{k}=v^{k}-w^{k}$ as a descent direction at $w^{k}$.
Take a suitable stepsize $\lambda_{k} \in (0,1]$,  set $w^{k+1}=w^{k}+\lambda_{k}p^{k}$
and $k=k+1$.
\medskip

The  stepsize can be found either with some one-dimensional minimization procedure
as in \cite{MF81} or with an inexact Armijo type linesearch; see also \cite{Pat98,BLM09}
for substantiation and further development.

The usefulness of this approach becomes clear if problem
(\ref{eq:3.1}) is (partially) decomposable, which is typical for very
large dimensional problems.
For instance, let
$$
\eta(w) = \sum \limits_{s \in N} \eta_{s}(w_{(s)}) \mbox{ and } \ W = \prod \limits_{s \in N} W_{s},
$$
where $w_{(s)} \in W_{s} \subset \mathbb{R}^{m_{s}}$, so that $m=\sum \limits_{s \in N} m_{s}$,
i.e., there is some concordant partition of the initial space $\mathbb{R}^{m}$. Then we have the problem
\begin{equation} \label{eq:3.3}
\min_{w \in \prod \limits_{s \in N} W_{s}} \to \left\{\mu(w) + \sum \limits_{s \in N} \eta_{s}(w_{(s)}) \right\},
\end{equation}
and (\ref{eq:3.2}) becomes equivalent to several
independent problems of the form
\begin{equation} \label{eq:3.4}
 \min\limits _{v_{(s)} \in W_{s}}  \to  \left\{
\left\langle v_{(s)}, \frac{\partial \mu(w^{k})}{\partial w_{(s)}}\right\rangle +\eta_{s}(v_{(s)})\right\};
\end{equation}
for $s \in N$.
The above descent method admits various component-wise iterative schemes; see e.g. \cite{Pat99}.

Our market equilibrium problem from the previous section written as VI (\ref{eq:2.3})
is reduced to problem
(\ref{eq:3.1}) in the case where the price functions are integrable, i.e.
$$
g_{is}(w)=\frac{\partial \mu(w)}{\partial x_{i}}, \ i \in I_{s}, \ \mbox{and} \
h_{js}(w)=-\frac{\partial \eta_{s}(w_{(s)})}{\partial y_{j}}, \ j \in J_{s}; \ s \in N.
$$
This is the case if these functions are separable, i.e. $g_{is}(w)=g_{is}(x_{is})$ for each $i \in I_{s}$ and
$h_{js}(w)=h_{js}(y_{js})$ for each $j \in J_{s}$, for all $s \in N$.
More precisely, VI (\ref{eq:2.3}) becomes the necessary optimality condition for
(\ref{eq:3.1}). The reverse assertion is true if the functions $\mu$ and $\eta$
are convex.

We now describe an adaptive cyclic component-wise PL method for problem
(\ref{eq:3.3}), which is some implementation of that from \cite{Kon16}.
For each point $w\in W$ and each $s \in N$, we define by
$V_{s}(w)$ the solution set of
the optimization problem:
$$
\min\limits _{v_{(s)} \in W_{s}}  \to  \left\{
\left\langle v_{(s)}, \frac{\partial \mu(w)}{\partial w_{(s)}}\right\rangle +\eta_{s}(v_{(s)})\right\};
$$
cf. (\ref{eq:3.4}). As above we suppose that the functions $\mu$ and $\eta$
are convex,   $\mu $ is smooth,
the set $W \subset \mathbb{R}^{m}$ is non-empty, convex, and compact.
Under these assumptions $V_{s}(w)$  is also non-empty, convex, and compact.
We define the gap function
$$
\varphi _{s} (w) = \max\limits_{v_{(s)} \in W_{s}}
\left\{ \left\langle w_{(s)}-v_{(s)}, \frac{\partial \mu(w)}{\partial w_{(s)}}\right\rangle +\eta_{s}(w_{(s)})-\eta_{s}(v_{(s)})\right\}
$$
for each $s \in N$. For brevity, set $f(w)=\mu(w)+\eta(w)$ and denote by $\mathbb{Z}_{+}$ the set of non-negative integers.
The optimal value of the function $f$ in (\ref{eq:3.1}) (or (\ref{eq:3.3})) will be denoted by $f^{*}$.
The adaptive cyclic PL method is described as follows.

\medskip \noindent
 {\bf Method (CPL).} \\ {\em Initialization:}
 Choose a point $z^{0}\in W$, numbers $\beta \in (0,1)$, $\theta
\in (0,1)$, and a sequence $\{\delta _{l}\} \searrow 0$. Set $l=1$.\\
{\em Step 0:} Set  $k=0$, $d=0$, $s=1$, $w^{0}=z^{l-1}$.\\
{\em Step 1:}  Solve problem
(\ref{eq:3.4}), find $v_{(s)}\in V_{s}(w^{k})$ and calculate
$\varphi _{s} (w^{k})$. If  $\varphi_{s}(w^{k}) \geq \delta _{l}$, take
$$
p_{(i)}^{k}= \left\{ {
\begin{array}{ll}
\displaystyle
v_{(s)}-w_{(s)}^{k} \quad & \mbox{if} \ i=s, \\
\mathbf{0} \quad & \mbox{if} \ i\neq s; \\
\end{array}
} \right.
$$
and go to Step 4.\\
{\em Step 2:}  Set $d=d+1$. If $d=n$, set $z^{l}=w^{k}$, $l=l+1$ and go to Step 0.
{\em (Restart)} \\
{\em Step 3:}  If $s=n$, set $s=1$, otherwise $s=s+1$. Afterwards  go to Step 1. \\
{\em Step 4:}  Determine $j$ as the smallest number in
$\mathbb{Z}_{+}$ such that
$$
 f(w^{k}+\theta ^{j} p^{k})
 \leq f(w^{k})-\beta \theta ^{j}\varphi_{s}(w^{k}),
$$
set $\lambda_{k}=\theta ^{j}$,
$w^{k+1}=w^{k}+\lambda_{k}p^{k}$,
$k=k+1$. If $s=n$, set $s=1$, otherwise $s=s+1$. Afterwards  go to Step 1.\\
\medskip

Thus, the method has two levels.
Each its outer iteration $l$ contains some number of inner iterations in
$k$ with the sequential verification of descent value for each component
with the fixed tolerance $\delta _{l}$. Completing each stage, which is marked as restart, leads to
decreasing the tolerance value. The basic properties of CPL are deduced along the same lines as in \cite{Kon16}.


\begin{proposition} \label{pro:3.1} Suppose in addition that the gradient
map the function $\mu$ is uniformly continuous on $W$. Then

 (i) the number of inner iterations at each outer iteration $l$  is finite;

 (ii) the sequence $\{z^{l}\}$ generated by Method (CPL)
 has limit points, all these limit points are solutions of
problem (\ref{eq:3.3}), besides,
$$
 \lim \limits_{l\rightarrow \infty} f (z^{l})=f^{*}.
$$
\end{proposition}

The line-search procedure in the method admits various modifications.
For instance,  we can take the exact one-dimensional minimization rule
instead of the current Armijo rule. If the gradient of the function $\mu$
is Lipschitz continuous, we can take fixed stepsize values and
remove the line-search procedure at all; see
\cite{Kon16} for more details.


\begin{remark} \label{rmk:3.1} Due to the presence of the control sequence
$\{\delta _{l}\}$, CPL differs essentially from the usual decomposition methods;
see e.g. \cite{Pat99,Mig04}. At the same time, this technique is rather usual for
for non-differentiable optimization methods; see e.g. \cite{BW75}. It was also applied in
iterative methods for linear inequalities \cite{Mac77} and
for decomposable variational inequalities \cite{Kon02a}.
\end{remark}


\section{A generalization of network equilibrium problems with elastic demands} \label{sc:4}

We now consider network flow equilibrium problems
with elastic (inverse) demands, which find various applications; see
\cite{Daf82}, \cite[Chapter IV]{Nag99} and references therein.

Let us be given a graph with finite sets of nodes $\mathcal{M}$ and
oriented arcs $\mathcal{A}$ which join the nodes so that any
arc $a=(i, j)$ has origin $i$ and destination $j$. Next,
among all the pairs of nodes of the graph we extract a subset of
origin-destination (O/D) pairs $\mathcal{N}$ of the form $s=(i \to
j)$. Each pair $s \in \mathcal{N}$ is associated  with the set of paths
$\mathcal{P}_{s}$ which connect the origin and destination for this
pair. Also, denote by $x_{p}$ the path flow for the path $p$.
Given a flow vector $ x =   (x_{p})_{ p \in \mathcal{P}_{s}, \, s \in \mathcal{N}
}$, one can determine the value of the arc
flow
\begin{equation} \label{eq:4.1}
f_{a}=     \sum _{s \in \mathcal{N}} \sum _{p  \in \mathcal{P}_{s}}
\alpha _{p a} x_{p}
\end{equation}
 for each arc  $a \in \mathcal{A}$, where
\begin{equation} \label{eq:4.2}
\alpha _{p a}  =\cases {
    1         &\textrm{ if arc $a$ belongs to path $p$,} \cr
    0   &\textrm{  otherwise.} \cr
          }
\end{equation}
If the vector $ f =   (f_{a})_{a \in \mathcal{A}}$ of arc flows is
known, one can determine the dis-utility value $c_{a}(f)$ for each
arc. Then one can compute the dis-utility value for each path
$p$:
\begin{equation} \label{eq:4.3}
g_{p}(x) =  \sum _{a  \in \mathcal{A}}   \alpha _{p a}  c_{a}(f).
\end{equation}
In the known elastic demand models, each (O/D) pair $s \in \mathcal{N}$ is
associated with one variable value of flow demand  
and hence one inverse demand (dis-utility) function; see
e.g. \cite[Chapter IV]{Nag99} and references therein.  
However, many active agents (users) with different dis-utility functions
may have the same physical location for many networks arising in applications.
For this reason, we now consider the generalization, where
each (O/D) pair $s \in \mathcal{N}$ may have several pairs of active users hence it
is associated  with the set of such pairs
$\mathcal{B}_{s}$ so that each pair of users $j \in \mathcal{B}_{s}$ has its particular
flow demand $y_{j}$ and dis-utility function $h_{j}$, which can be in principle dependent of the flow demand
$y$, i.e. $y = (y_{j})_{j \in \mathcal{B}_{s}, s \in \mathcal{N}}$.
 Then one can define the feasible set of flows:
\begin{equation} \label{eq:4.4}
W = \left\{ w=(x,y)
  \ \vrule \ \begin{array}{c}
   \sum _{p  \in \mathcal{P}_{s}} x_{p} =\sum _{j  \in \mathcal{B}_{s}} y_{j}, \ x_{p} \geq 0, \
     p \in \mathcal{P}_{s},  \\
   0 \leq y_{j} \leq \gamma_{j}, \ j \in \mathcal{B}_{s}; \ s \in \mathcal{N}
  \end{array}
  \right\}.
\end{equation}
We say that
a feasible flow / demand pair $ (x^{*}, y^{*} ) \in W$ is an {\em
equilibrium point} if it satisfies the following conditions:
\begin{equation} \label{eq:4.5}
\forall s \in \mathcal{N}, \ \exists \lambda_{s} \ \mbox{such that}
\ g_{p}(x^{*}) \left\{
\begin{array}{ll}
\geq \lambda_{s}  \quad & \mbox{if} \quad x^{*}_{p} = 0,  \\
=\lambda_{s}   \quad & \mbox{if} \quad x^{*}_{p}  >0,
\end{array}
\right. \quad   \forall p \in \mathcal{P}_{s};
\end{equation}
and
\begin{equation} \label{eq:4.6}
 h_{j}(y^{*}) \left\{
\begin{array}{ll}
\leq \lambda_{s}  \quad & \mbox{if} \quad y^{*}_{j} = 0,  \\
=\lambda_{s}   \quad & \mbox{if} \quad y^{*}_{j} \in (0,\gamma_{j}), \\
\geq \lambda_{s}
 \quad & \mbox{if} \quad y^{*}_{j} =\gamma_{j},
\end{array}
\right. \quad   \forall j \in \mathcal{B}_{s}.
\end{equation}
Clearly, the equilibrium conditions in
(\ref{eq:4.5})--(\ref{eq:4.6}) represent some implementation of
the multi-commodity two-sided  market equilibrium model (\ref{eq:2.1})--(\ref{eq:2.2}),
where each commodity is associated with an (O/D) pair $s \in \mathcal{N}$,
its set of traders (carriers) with price functions $g_{p}(x)$ is represented by the paths $p \in
\mathcal{P}_{s}$, whereas
its set of buyers with price functions $h_{j}(y)$ is represented by the pairs of users $j \in
\mathcal{B}_{s}$.  We observe that the prices here
are not fixed, the dependence of volumes for offer price
functions $g_{p}$ is given in (\ref{eq:4.1})--(\ref{eq:4.3}) and
caused by the complexity of the system topology and by the
fact that carriers of different (O/D) pairs can utilize the same
links (arcs).

We now show that conditions (\ref{eq:4.4})--(\ref{eq:4.6}) can be
equivalently rewritten in the form of a VI: Find a pair $
(x^{*},y^{*}) \in W$ such that
\begin{equation} \label{eq:4.7}
\sum _{s \in \mathcal{N}} \sum _{p  \in \mathcal{P}_{s}}
g_{p}(x^{*}) (x_{p}- x^{*}_{p})-\sum _{s \in \mathcal{N}} \sum _{j  \in \mathcal{B}_{s}}
h_{j}(y^{*}) (y_{j}- y^{*}_{j})  \geq 0       \quad \forall (x,y)
\in W.
\end{equation}

\begin{proposition} \label{pro:4.1}
 A pair $(x^{*},y^{*})\in W$ solves VI (\ref{eq:4.7}) if and only if it
satisfies conditions (\ref{eq:4.5})--(\ref{eq:4.6}).
\end{proposition}
{\bf Proof.}
Writing the usual necessary and sufficient optimality conditions
(see \cite[Proposition 11.7]{Kon07a}) for problem (\ref{eq:4.7}), we
obtain that there exist $x^{*} \geq \mathbf{0}$, $y^{*} \in [\mathbf{0},
\gamma]$, and $\lambda$ such that
\begin{eqnarray*}
\displaystyle
 & & \sum \limits_{p  \in \mathcal{P}_{s}} (g_{p}(x^{*})-\lambda_{s})(x_{p}-x_{p}^{*})\geq 0
       \quad \forall x_{p}\geq 0, \ p  \in \mathcal{P}_{s}, \ s \in \mathcal{N};   \\
 & &  \sum _{k  \in \mathcal{B}_{s}}\left(\lambda_{s}-h_{j}(y^{*}) \right)(y_{j}- y^{*}_{j})\geq 0
                       \quad  \forall y_{j} \in (0,\gamma_{j}), \ s \in \mathcal{N}; \\
                       & &  \sum _{p  \in \mathcal{P}_{s}} x^{*}_{p} =\sum _{k  \in \mathcal{B}_{s}} y^{*}_{j}, \quad j \in \mathcal{B}_{s}, \ s \in \mathcal{N};
\end{eqnarray*}
where $\lambda = (\lambda_{s})_{s \in \mathcal{N}}$. However, the
first and second relations are clearly equivalent to
(\ref{eq:4.5})--(\ref{eq:4.6}).
\QED

If each (O/D) pair is attributed to only one pair of users, we obtain the
custom network  equilibrium problems with elastic (inverse) demands; see
e.g. \cite[Chapter IV]{Nag99}. If all the (O/D) traffic demands in this model are not
restricted with upper bounds, we obtain the model considered in \cite{Daf82}. Let us insert the same condition in our model:
\begin{equation} \label{eq:4.8}
\gamma_{j}=+\infty \quad \forall j \in \mathcal{B}_{s}, \ s \in \mathcal{N}.
\end{equation}
 Then (\ref{eq:4.6}) reduces to the following condition:
\begin{equation} \label{eq:4.9}
 h_{j}(y^{*}) \left\{
\begin{array}{ll}
\leq \lambda_{s}  \quad & \mbox{if} \quad y^{*}_{j} = 0,  \\
=\lambda_{s}   \quad & \mbox{if} \quad y^{*}_{j} >0;
\end{array}
\right. \quad    \forall j \in \mathcal{B}_{s}.
\end{equation}

We can also write some other equivalent network
equilibrium conditions, for instance,
\begin{equation} \label{eq:4.10}
\begin{array}{c}   g_{p}(x^{*})-h_{j}(y^{*})
 \left \{ \begin{array}{lllcc}
 = 0  & \ \mbox{if} &  x_{p}^{*} > 0 & \mbox{ and } & y^{*}_{j} > 0,\\
 \geq 0  & \  \mbox{if} &  x_{p}^{*} = 0  & \mbox{ or } & y^{*}_{j} = 0;
 \end{array} \right. \\
 \quad \forall p \in \mathcal{P}_{s}, \ j \in \mathcal{B}_{s}, \ s \in \mathcal{N}.
 \end{array}
\end{equation}

\begin{proposition} \label{pro:4.2}  Let (\ref{eq:4.8}) hold.
Then, for any pair $(x^{*},y^{*})\in W$, condition (\ref{eq:4.10})
is equivalent to (\ref{eq:4.5})  and (\ref{eq:4.9}).
\end{proposition}
{\bf Proof.}
Take an arbitrary pair $s \in \mathcal{N}$. Suppose a pair
$(x^{*},y^{*})\in W$ satisfies conditions (\ref{eq:4.5})  and (\ref{eq:4.9}).
 Then, for any $p \in \mathcal{P}_{s}$ and $j \in \mathcal{B}_{s}$, the relations
$x^{*}_{p} >0$ and $y^{*}_{j} > 0$ imply $g_{p}(x^{*})=\lambda_{s}=h_{j}(y^{*})$. Next,
each of the relations $x^{*}_{p} =0$ or $y^{*}_{j} = 0$ implies
$g_{p}(x^{*})\geq \lambda_{s}\geq h_{j}(y^{*})$, and
(\ref{eq:4.10}) holds true.

Conversely, suppose a pair $(x^{*},y^{*})\in W$ satisfies conditions
(\ref{eq:4.10}). Fix any $s \in \mathcal{N}$ and set
$$
 \alpha'=\min_{p \in \mathcal{P}_{s}} g_{p}(x^{*}), \ \alpha''=\max_{j \in \mathcal{B}_{s}}
 h_{j}(y^{*}),
$$
then $\alpha' \geq \alpha''$.  If $x_{p}^{*}=0$ for all $p \in \mathcal{P}_{s}$,
then $y^{*}_{j}=0$ for all $j \in \mathcal{B}_{s}$ and conversely. Then taking any
$\lambda_{s} \in [\alpha'', \alpha']$ yields (\ref{eq:4.5})  and (\ref{eq:4.9}).
Otherwise,
there exists at least one pair of indices $p \in \mathcal{P}_{s}$, $j \in \mathcal{B}_{s}$ such
that $x_{p}^{*}>0$ and $y^{*}_{j}>0$. Then setting $\lambda_{s} =\alpha'
= \alpha''$ again yields (\ref{eq:4.5})  and (\ref{eq:4.9}).\QED

It is easy to see that conditions
(\ref{eq:4.10}) can be replaced with the following:
\begin{equation} \label{eq:4.11}
\begin{array}{c}   g_{p}(x^{*})-h_{j}(y^{*})
 \left \{ \begin{array}{lllcc}
 > 0  & \ \Longrightarrow &  x_{p}^{*} = 0  & \mbox{ or } & y^{*}_{j} = 0, \\
 \geq 0  & \  \Longleftrightarrow  &   x_{p}^{*} \geq 0 & \mbox{ and } & y^{*}_{j} \geq 0;
 \end{array} \right. \\
 \quad \forall p \in \mathcal{P}_{s}, \ j \in \mathcal{B}_{s}, \ s \in \mathcal{N}.
 \end{array}
\end{equation}

\begin{proposition} \label{pro:4.3}
Let (\ref{eq:4.8}) hold.
Then, for any pair $(x^{*},y^{*})\in W$, condition (\ref{eq:4.11})
is equivalent to (\ref{eq:4.5})  and (\ref{eq:4.9}).
\end{proposition}

The equivalent VI formulation of network equilibrium problems enables us to
obtain the existence of solutions rather easily.
The feasible set $W$ of
the network equilibrium problem defined in
(\ref{eq:4.4}) is bounded if $\gamma_{j}<+\infty$ for all $j \in \mathcal{B}_{s}$, $s \in \mathcal{N}$.
Then VI (\ref{eq:4.7}) and hence the
equivalent network equilibrium problem are
solvable if all the mappings $c_{a}$, $a \in \mathcal{A}$ and $h_{j}$, $j \in \mathcal{B}_{s}$, $s \in \mathcal{N}$ are
continuous. Let us turn to the above pure unbounded case (\ref{eq:4.8}).
Then the feasible set $W$ in
(\ref{eq:4.4}) is unbounded.  We now deduce a new
existence result for VI (\ref{eq:4.7}) and hence for the equivalent
network equilibrium problem by a direct application of Proposition \ref{pro:2.2}.
We need the proper following coercivity condition; cf. {\bf (C)}.

{\bf (C1)} {\em There exists a number $r>0$ such that for any point
$w=(x,y) \in W$ and for each $s \in \mathcal{N}$ it holds that}
$$
 \exists j \in \mathcal{B}_{s}, \ y_{j} >r \Longrightarrow \exists p \in \mathcal{P}_{s}  \  \mbox{such that} \
   x_{p}>0 \ \mbox{and} \  g_{p}(x)\geq h_{j}(y).
$$

We observe that condition  {\bf (C1)}  implies condition {\bf (C)} for VI (\ref{eq:4.7})
and we obtain the desired  existence result.


\begin{theorem} \label{thm:4.1} Suppose that (\ref{eq:4.8}) holds, the set $W$ defined
in (\ref{eq:4.4}) is nonempty, all the functions $c_{a}$ and $h_{j}$ are continuous
for all $a \in \mathcal{A}$, $j \in \mathcal{B}_{s}$, and $s \in \mathcal{N}$. If condition {\bf (C1)}
is fulfilled, then VI (\ref{eq:4.7}) has a solution.
\end{theorem}


\section{Implementation of partial linearization methods for integrable network equilibrium problems} \label{sc:5}

In Section \ref{sc:3}, several versions of partial linearization (PL) methods
for special decomposable optimization  problems over Cartesian product sets
were described for the general multi-commodity market equilibrium model
of Section \ref{sc:2} in the integrable case. Hence, PL methods
can be also applied to integrable network equilibrium problems
with elastic demands of Section \ref{sc:4}.

Therefore, we now will suppose that all the functions $c_{a}$ and $h_{j}$ are continuous and separable,
i.e., $c_{a}(f)=c_{a}(f_{a})$ and $h_{j}(y)=h_{j}(y_{j})$. Besides, we assume that
$c_{a}(f_{a})$ and $-h_{j}(y_{j})$ are monotone
increasing functions. Next, we assume that
$$
\gamma_{j}<+\infty \quad \forall j \in \mathcal{B}_{s}, \ s \in \mathcal{N};
$$
then the feasible set $W$ is  non-empty, convex, and compact and
\begin{eqnarray*}
   W &=& \prod_{s \in \mathcal{N}} W_{s}, \ \mbox{where} \\
   W_{s}&=& \left\{ w_{(s)}=(x_{(s)},y_{(s)}) \ \vrule \
\begin{array}{l}
 \sum _{p  \in \mathcal{P}_{s}} x_{p} =\sum _{j  \in \mathcal{B}_{s}} y_{j}, \\
 x_{p} \geq 0, p \in \mathcal{P}_{s}, 0 \leq y_{j} \leq \gamma_{j}, \ j \in \mathcal{B}_{s}
\end{array}
\right\}; \\
 && \mbox{for} \ s \in \mathcal{N}.
\end{eqnarray*}
Here $x_{(s)}=(x_{p})_{p  \in \mathcal{P}_{s}}$, $y_{(s)}=(y_{j})_{j \in \mathcal{B}_{s}}$.

Due to the separability of the functions $c_{a}$ and $h_{j}$, their continuity
implies integrability, i.e., then there exist functions
$$
\mu_a(f_a) = \int \limits_{0}^{f_a} c_a(t) dt \ \forall a \in \mathcal{A}, \
\eta_j(y_{j}) = \int \limits_{0}^{v_j} h_j(t) dt \  \forall j  \in \mathcal{B}_{s}, \ s \in \mathcal{N}.
$$
Taking into account (\ref{eq:4.1}), we see that VI (\ref{eq:4.7}) gives a necessary and sufficient
 optimality condition for the following optimization problem:
\begin{equation}\label{eq:5.1}
\min \limits_{(x, y) \in W} \rightarrow
\left\{ \sum\limits_{a \in \mathcal{A}} \mu_a(f_a) -
 \sum\limits_{s \in \mathcal{N}} \sum\limits_{j  \in \mathcal{B}_{s}} \eta_j(y_{j})\right\}.
\end{equation}
 However, this problem falls into the basic format
(\ref{eq:3.3}) and the suggested PL methods can be applied to (\ref{eq:5.1}).

We describe the solution of the basic direction finding problem (\ref{eq:3.4}).
It now consists in finding an element
$\bar w_{(s)}=(\bar x_{(s)},\bar y_{(s)}) \in W_{s}$,
which solves the optimization problem
\begin{equation}\label{eq:5.2}
\min \limits_{(x_{(s)}, y_{(s)}) \in W_{s}} \rightarrow
\left\{ \sum _{p  \in \mathcal{P}_{s}}
g_{p}(x^{k})x_{p} - \sum _{j  \in \mathcal{B}_{s}}
  \eta_j(y_{j})\right\}
\end{equation}
for some selected pair $s \in \mathcal{N}$. The solution of (\ref{eq:5.2}) can be found with the simple procedure below,
which is based on optimality conditions (\ref{eq:4.5})--(\ref{eq:4.6}).

First we calculate the shortest path $q \in \mathcal{P}_{s}$  with the minimal
cost. Set $\tilde \lambda_{s}=g_q(x^k)$, $\bar x_{p}=0$ for all $p \in \mathcal{P}_{s}$.

For each $j  \in \mathcal{B}_{s}$ we verify three possible cases.

{\em Case 1.} If $h_{j}(0) \le \tilde \lambda$, then set $\bar y_{j}=0$.
Otherwise go to Case 2.

{\em Case 2.} If $h_{j}(\gamma_{j}) \ge \tilde \lambda$,
set $\bar y_{j}=\gamma_{j}$, $\bar x_{q}=\bar x_{q}+\gamma_{j}$.
Otherwise go to Case 3.

{\em Case 3.} We have $h_{j}(\gamma_{j}) < \tilde \lambda < h_{j}(0)$. By continuity of
$h_{j}$, we find the value $\bar y_{j} \in [0, \gamma_{j}]$  such that
$h_{j}(\bar y_{j}) =\tilde \lambda$, set
$\bar x_{q}=\bar x_{q}+\bar y_{j}$.

Therefore, the suggested PL methods can be implemented rather easily.


\section{Application of market models to resource allocation in wireless networks} \label{sc:6}

In contemporary wireless networks, increasing demand of services
leads to serious congestion effects, whereas significant network
resources (say, bandwidth and batteries capacity) are utilized
inefficiently for systems with fixed allocation rules.
This situation forces one to apply more
flexible market type allocation mechanisms.
Due to the presence of conflict of interests, most papers on
allocation mechanisms are devoted to pure game-theoretic models
reflecting imperfect competition; see, e.g., \cite{LZ09,RAR10}.
However,  certain lack of information
about the participants is typical for wireless
telecommunication networks (see, e.g., \cite{IK10,RAR10}), and some other
market models may be suitable here because they can be
utilized under minimal information requirements on involved users.

We now consider the problem of allocation of services of several
competitive wireless network providers for a large number of users,
which is very essential for contemporary communication systems.
This problem was investigated in
\cite{HTW07,KIAK12,MTV12,ZZ13} for wired and wireless network
settings, where game-theoretic models for competitive providers were
presented. An alternative model, which is based on some VI formulation and uses
proper equilibrium conditions, was suggested for this problem in
\cite[Section 6]{Kon15c}. We now propose its extension that admits different
kinds of users' behavior.

Namely, we suppose that there are $m$
wireless network providers and that all the users are divided into
$n$ classes, that is, the users belonging to the same class $j$ are
considered as one service consumer with a price function
$h_{j}(y_{j})$ and a scalar bid volume $y_{j} \in [0, \beta_{j}]$
for $j \in N=\{1, \ldots, n\}$. Next, each provider $i$ announces
his/her price function $b_{i}(x_{i})$ depending on the offer volume
$x_{i} \in [0, \alpha _{i}]$ for $i \in M=\{1, \ldots, m\}$.
However, such joint consumption of wireless network resources
yields the additional dis-utility $l_{i}(x)$ for users consuming
resources of provider $i$, where $x=(x_{1}, \ldots, x_{s})^{\top}$;
see \cite{KIAK12,MTV12,ZZ13} for more detail. Hence, the actual price function of
provider $i$ for users becomes $g_{i}(x)=b_{i}(x_{i})+l_{i}(x)$. We
can thus define the feasible set of offer/bid values
$$
\displaystyle D=\left\{ (x, y) \ \vrule \
\begin{array}{l} x_{i} \in [0, \alpha _{i}], i \in M, \\ y_{j} \in
[0, \beta_{j}], j \in N;
\end{array}
\sum \limits_{i \in M} x_{i} =  \sum \limits_{j \in N}
y_{j};\right\};
$$
where $y=(y_{1}, \ldots, y_{n})^{\top}$. Then we can write the two-sided equilibrium
problem that consists in finding a
feasible pair $(\bar x, \bar y) \in D$ and a price $\lambda$ such
that
\begin{equation} \label{eq:6.1}
\begin{array}{l}
g_{i}(\bar x) \left\{
\begin{array}{l}
\geq \lambda, \  \mbox{if} \ \bar x_{i} = 0,  \\
=\lambda,  \  \mbox{if} \ \bar x_{i} \in (0, \alpha_{i}), \\
\leq \lambda, \  \mbox{if} \ \bar x_{i}=\alpha _{i},
\end{array}
\right. \\
 i \in M;
\end{array}
\begin{array}{l}
h_{j}(\bar y_{j}) \left\{
\begin{array}{l}
\leq \lambda, \  \mbox{if} \ \bar y_{j}= 0, \\
=\lambda,  \  \mbox{if} \ \bar y_{j} \in (0, \beta_{j}), \\
\geq \lambda, \  \mbox{if} \ \bar y_{j}= \beta_{j},
\end{array}
 \right. \\
 \ j \in N.
\end{array}
\end{equation}
Clearly, it is a particular case of those in (\ref{eq:2.1})--(\ref{eq:2.2}).
Due to Proposition \ref{pro:2.1},  (\ref{eq:6.1}) can be replaced with
the equivalent VI: Find  $(\bar x, \bar y) \in D$ such that
\begin{equation} \label{eq:6.2}
\sum \limits_{i \in M} g_{i} (\bar x) (x_{i} - \bar x_{i}) - \sum
\limits_{j \in N} h_{j} (\bar y_{j}) (y_{j} - \bar y_{j}) \geq 0
\quad \forall (x, y) \in D.
\end{equation}
This property enables us to establish existence of solutions for the
above problem and develop efficient iterative solution methods.
In fact, if all the price functions are continuous and
the set $D$ is nonempty and bounded, then VI (\ref{eq:6.2}) has a solution. In
the unbounded case, some coercivity condition is necessary. For
instance, let us consider the case where $\alpha _{i}=+\infty$ for
$i \in M$ and $\beta_{j}=+\infty$ for $j \in N$ and take the
following condition; cf. {\bf (C)}.

{\bf (C2)} {\em There exists a number $r>0$ such that for any pair
$(x,y) \in D$ it holds that}
$$
y_{l} >r \Longrightarrow \exists k \in M \  \mbox{such that} \  x_{k}>0 \ \mbox{and} \
 g_{k}(x)\geq h_{l}(y_{l}).
$$

Clearly, {\bf (C2)}  implies {\bf (C)} for VI (\ref{eq:6.2})
 and  Proposition \ref{pro:2.2} provides the existence
result.


\begin{theorem} \label{thm:6.1} Suppose that the set $D$ is nonempty,
the functions $g_{i}$ and $h_{j}$ are continuous for all $i \in M$, $j
\in N$. If condition {\bf (C2)} is fulfilled, then VI (\ref{eq:6.2}) has a
solution.
\end{theorem}


\section{The partial linearization method
for resource allocation problems in wireless networks} \label{sc:7}

Iterative solution methods for solving  VI of form (\ref{eq:6.2}) in general
require additional monotonicity assumptions for convergence; see e.g. \cite{Pat99,Kon07a,Kon13a}.
Additional solution methods appear in the  integrable case where
$$
g_{i}(x)=\frac{\partial \mu(x)}{\partial x_{i}}, \ i \in M;  \
h_{j}(y_{j})=-\eta'_{j}(y_{j}), \ j \in N.
$$
Then, VI (\ref{eq:6.2}) gives the optimality condition for the optimization
problem:
\begin{eqnarray}
& & \min \limits _{w \in D} \to f(w),  \label{eq:7.1}\\
  f(w)=f(x, y) &=&  \left\{\mu(x)+\eta(y)\right\},  \ \eta(y)=\sum \limits_{j \in N} \eta_{j} (y_{j}); \nonumber
\end{eqnarray}
cf. (\ref{eq:3.1}) and (\ref{eq:3.3}).
In particular, conditional gradient, gradient projection, and Uzawa
type methods then can be utilized; see e.g. \cite{Kon13,Kon15e}.
We now only describe a way to implement the custom PL method since the problem is not separable.
We suppose in addition that the function $\mu$ is smooth and convex, $\alpha _{i}=+\infty$ for all
$i \in M$, and $0 \leq \beta_{j}<+\infty$ for all $j \in N$.
Then the feasible set $D$ is  non-empty, convex, and compact.

For more clarity, we rewrite the PL method for problem (\ref{eq:7.1}).
We define the gap function
$$
\varphi (w) = \varphi (x, y) = \max\limits_{(x',y') \in D}
\left\{ \langle x-x', \mu'(x)\rangle +\eta(y)-\eta(y')\right\}.
$$

\medskip \noindent
 {\bf Method (PL).} \\
Choose a point $w^{0}\in D$, numbers $\beta \in (0,1)$ and $\theta
\in (0,1)$, set $k=0$. At the $k$-th iteration, $k=0,1,\ldots$, we have
a point $w^{k}\in D$. Find a solution $v^{k}=(\bar x^{k},\bar y^{k})$ of the problem
\begin{equation} \label{eq:7.2}
\min_{v \in D} \to \left\{\langle \mu'(x^{k}),v\rangle + \eta(v)\right\}.
\end{equation}
If $v^{k}=w^{k}$, stop.
Otherwise set $d^{k}=v^{k}-w^{k}$, find $p$ as the smallest number in
$\mathbb{Z}_{+}$ such that
$$
 f(w^{k}+\theta ^{p} d^{k})
 \leq f(w^{k})-\beta \theta ^{p} \varphi (w^{k}),
$$
set $\sigma_{k}=\theta ^{p}$,
$w^{k+1}=w^{k}+\sigma_{k}d^{k}$, and $k=k+1$.
\medskip

The solution of the basic direction finding problem (\ref{eq:7.2})
 can also be found with the simple procedure,
which is similar to that from Section \ref{sc:5} and based on the optimality
conditions.

First we calculate an index $q \in M$ that corresponds to the minimal value
$$
 g_{q}(x^{k})=\min_{i \in M}g_{i}(x^{k})
$$
and set $\tilde \lambda=g_q(x^k)$, $\bar x^{k}_{i}=0$ for all $i \in M$.

For each $j  \in N$ we verify three possible cases.

{\em Case 1.} If $h_{j}(0) \le \tilde \lambda$, then set $\bar y^{k}_{j}=0$.
Otherwise go to Case 2.

{\em Case 2.} If $h_{j}(\beta_{j}) \ge \tilde \lambda$,
set $\bar y^{k}_{j}=\beta_{j}$, $\bar x^{k}_{q}= \bar x^{k}_{q}+\beta_{j}$.
Otherwise go to Case 3.

{\em Case 3.} We have $h_{j}(\beta_{j}) < \tilde \lambda < h_{j}(0)$. By continuity of
$h_{j}$, we find the value $\bar y^{k}_{j} \in [0, \beta_{j}]$  such that
$h_{j}(\bar y^{k}_{j}) =\tilde \lambda$, set
$\bar x^{k}_{q}= \bar x^{k}_{q}+\bar y^{k}_{j}$.

Let us now consider the case where $0 \leq \alpha _{i}<+\infty$ for all
$i \in M$ and $0 \leq \beta_{j}<+\infty$ for all $j \in N$.
Then the feasible set $D$ is also non-empty, convex, and compact.
Hence, the above PL method can be applied to (\ref{eq:7.1}), however,
we should then take more complex procedures for solution of
 problem (\ref{eq:7.2}). However, we can eliminate the upper bounds
for the variables $x_{i}$ via a suitable penalty approach.

For instance, replace problem (\ref{eq:7.1}) with
the sequence of auxiliary problems of the form
\begin{eqnarray}
& & \min \limits _{w \in D} \to \Phi(w,\tau),  \label{eq:7.3}\\
  \Phi(w,\tau) &=&  \mu(x)+\tau\varphi(x)+\eta(y),
  \ \varphi(x)=0.5 \sum \limits_{i \in M} \max\{x_{i}-\alpha_{i},0\}^{2}; \nonumber
\end{eqnarray}
where $\tau>0$ is a penalty parameter, the functions $\mu$ and $\eta$ are defined as above.
Under the standard assumptions the sequence of solutions of (\ref{eq:7.3}) will approximate a
solution of (\ref{eq:7.1}) if $\tau \to +\infty$; see e.g. \cite{Kon13a}.
Next, each problem (\ref{eq:7.3}) has the previous format without 
the upper bounds for the variables $x_{i}$. Hence, we can apply directly the above version of
the PL method to (\ref{eq:7.3}) with replacing $f(w)$ by $\Phi(w,\tau)$.
Clearly, (\ref{eq:7.2}) is replaced by
$$
\min_{v \in D} \to \left\{\langle \mu'(x^{k})+\tau\varphi'(x^{k}),v\rangle + \eta(v)\right\}.
$$
We also have to substitute each function $g_{i}(x)$
with $\tilde g_{i}(x)=g_{i}(x)+\tau \max\{x_{i}-\alpha_{i},0\}$ in the procedure
of finding its solution. This gives us an alternative way to solve
such resource allocation problems in wireless  networks.


\section{Computational experiments with network equilibrium test problems} \label{sc:8}

In order to compare the performance of the PL methods
 we carried out preliminary series of computational experiments
on network equilibrium test problems
of form (\ref{eq:4.4})--(\ref{eq:4.6}) or (\ref{eq:4.7}). We took their adjustment
described in Section \ref{sc:5}.

For comparison we took proper extensions of the known
test examples of network equilibrium problems
with elastic demands, namely, each (O/D) pair
was associated  with two pairs of active users.
We used the arc cost functions $c_{a}(f_{a})= 1 + f_a$
for all $a \in \mathcal{A}$ and the minimal path cost (dis-utility) functions
$h_{j1(s)}(y_{j1})= 30 - 0.5 y_{j1(s)}$ and $h_{j2(s)}(y_{j2(s)})= 28 - 0.3 y_{j2(s)}$,
 where $\mathcal{B}_{s}=\{j1(s),j2(s)\}$ for all $s \in \mathcal{N}$.
We took
$$
\Delta_{k}=\varphi(w^{k})=\sum \limits_{s \in \mathcal{N}} \varphi_{s}(w^{k})
$$
as accuracy measure for the methods.
Both the PL and CPL methods were implemented with the Armijo line-search rule
where  $\beta =\theta=0.5$. Due to the above description we see that we can take
the total number of blocks where the linesearch procedure was utilized
as unified complexity measure for both the methods,
which will be called block iterations. Hence we reported this value in the tables for
attaining different accuracies. The methods were
implemented in C++ with double precision arithmetic.

The topology of Example 1 was taken from \cite{BG82}.
The graph contains 25 nodes, 40 arcs, and 5 O/D pairs. We used two rules for
changing the parameter $\delta _{l}$ with $\delta _{0}=10$ in CPL.
The performance results are given in Table \ref{tbl:1}.
\begin{table}[h]
\caption{Example 1. The numbers of block iterations} \label{tbl:1}
\begin{center}
\begin{tabular}{|r|c|c|c|}
\hline \ accuracy & PL & CPL & CPL \\
\hline \  &  & $\delta _{l+1}=\delta _{l}/2$ & $\delta _{l}=\delta _{0}/l$ \\
\hline 0.2  &  4970 &  4427 &  3519 \\
\hline 0.1  & 10785 &  8747 &  6411 \\
\hline 0.05 & 21260 & 17284 & 13425 \\
\hline
\end{tabular}
\end{center}
\end{table}

The topology of Example 2 was taken from \cite[Network 26]{Nag84}.
The graph contains 22 nodes, 36 arcs, and 12 O/D pairs. We used the rule
$\delta _{l}=\delta _{0}/l$ with $\delta _{0}=10$ in CPL.
The performance results are given in Table \ref{tbl:2}.
\begin{table}[h]
\caption{Example 2. The numbers of block iterations} \label{tbl:2}
\begin{center}
\begin{tabular}{|r|c|c|}
\hline \ accuracy & PL & CPL \\
\hline 0.2  &  420 & 233 \\
\hline 0.1  & 468 & 246 \\
\hline 0.05 & 504 & 256 \\
\hline
\end{tabular}
\end{center}
\end{table}

In Example 3, the data were generated randomly.
The graph contained 20 nodes, 114 arcs, and 10 O/D pairs.
We used the rule $\delta _{l}=\delta _{0}/l$ with $\delta _{0}=10$ in CPL. The results are given in Table \ref{tbl:3}.
\begin{table}[h]
\caption{Example 3. The numbers of block iterations} \label{tbl:3}
\begin{center}
\begin{tabular}{|r|r|r|}
\hline \ accuracy & PL & CPL \\
\hline 1    &  135730 &  106308 \\
\hline 0.5  &  271830 &  217932 \\
\hline 0.2  &  662220 &  531032 \\
\hline 0.1  & 1329910 & 1082449 \\
\hline
\end{tabular}
\end{center}
\end{table}
In all the cases, CPL showed certain preference
over PL in the number of block iterations.


\section{Conclusions} \label{sc:9}

We considered the general market model with many divisible
commodities and price functions of participants and established
existence results for this problem under natural coercivity
conditions in the case of an unbounded feasible set.
We described extensions
of the known network flow equilibrium problems with elastic
demands and  a resource allocation problem in wireless communication
networks and showed they are particular cases of the presented market
model.  This property enabled us to obtain new existence results for all these models as some
adjustments of that for the general market model.
Besides, under certain integrability
conditions the market model can be reduced to an optimization problem.
We suggested a new cyclic version of the partial linearization (PL) method
for its decomposable case. We suggested ways for implementation of
the PL method to solve the network equilibrium  problems and
resource allocation problems in wireless communication
networks.

\section*{Acknowledgements}
This work was supported by the RFBR grant, project No. 16-01-00109a
and by grant No. 297689 from Academy of Finland. The author is grateful to Olga Pinyagina for
her assistance in carrying out computational experiments.


\end{document}